\documentclass[a4paper,fleqn]{cas-dc}

\usepackage[authoryear]{natbib}

\usepackage{graphicx}
\usepackage{xcolor}

\usepackage{amsmath, amssymb, mathtools, amsfonts}
\usepackage{algorithm}
\usepackage[noend]{algpseudocode}
\renewcommand{\Comment}[2][.15\linewidth]{%
  \leavevmode\hfill\makebox[#1][l]{\textcolor{gray}{// #2}}}
\usepackage[short,nocomma]{optidef}
\usepackage{booktabs}   
\usepackage{pifont}
\newcommand{\cmark}{\ding{51}}
\newcommand{\xmark}{\ding{55}}
\usepackage{microtype}
\usepackage{inconsolata}

\usepackage{amsthm}
\theoremstyle{plain}

\newtheorem{lemma}[]{Lemma}

\theoremstyle{definition}

\newtheorem{remark}[]{Remark}

\algdef{SE}[IF]{NoThenIf}{EndIf}[1]{\algorithmicif\ #1}{\algorithmicend\ \algorithmicif}%

\algblockdefx{MRepeat}{EndRepeat}{\textbf{repeat}}{}
\algnotext{EndRepeat}

\newcommand{\bmat}{\begin{bmatrix}}
\newcommand{\emat}{\end{bmatrix}}
\newcommand{\dpartial}[2]{\frac{\partial{#1}}{\partial{#2}}}

\newcommand{\dtotaltfrac}[2]{\tfrac{\mrm{d} #1}{\mrm{d} #2}}

\DeclareMathOperator{\blkdiag}{blkdiag}
\newcommand{\argmin}{\arg\min}
\DeclareMathOperator{\tr}{tr}
\DeclareMathOperator{\EXV}{\mathbb{E}}		
\DeclareMathOperator{\vecop}{vec}		

\DeclareMathSymbol{\shortminus}{\mathbin}{AMSa}{"39}

\newcommand{\R}{\mathbb{R}}
\newcommand{\SSS}{{\mathbb{S}}}

\newcommand{\calE}{\mathcal{E}}

\newcommand{\bigO}{\mathcal{O}}
\newcommand{\defeq}{:=}
\newcommand{\eye}{I}  

\newcommand{\mrm}{\mathrm}

\usepackage{accents}
\newlength{\dhatheight}

\renewcommand{\eqdef}{=:}

\begin{document}
\let\WriteBookmarks\relax
\def\floatpagepagefraction{1}
\def\textpagefraction{.001}

\shorttitle{Riccati-ZORO: Efficient online optimization of internal feedback laws in robust and stochastic MPC}
\shortauthors{F. Messerer et al.}

\title [mode = title]{%
Riccati-ZORO:
An efficient algorithm for heuristic online optimization of internal feedback laws in robust and stochastic model predictive control%
}  

\tnotemark[1] 
\tnotetext[1]{%
  This research was supported by DFG via projects 504452366 (SPP 2364), 560056112 (robust MPC), 535860958 (ALeSCo) and 525018088 (MAWERO), and by BMWK via 03EN3054B.}

\author[1]{Florian Messerer}%
[orcid=0000-0002-6984-9585]
\cormark[1]
\ead{florian.messerer@imtek.uni-freiburg.de}
\credit{Conceptualization, Methodology, Software, Validation, Formal analysis, Investigation, Data Curation, Writing -- Original Draft, Writing -- Review \& Editing, Visualization, Project Administration}
\author[1]{Yunfan Gao}%
[orcid=0009-0001-8100-0605]
\ead{yunfan.gao@imtek.uni-freiburg.de}
\credit{Software, Validation, Investigation, Writing -- Review \& Editing, Visualization}
\author[1]{Jonathan Frey}%
[orcid=0000-0003-2771-4209]
\ead{jonathan.frey@imtek.uni-freiburg.de}
\credit{Software, Data Curation, Writing -- Review \& Editing}
\author[1,2]{Moritz Diehl}%
[orcid=0000-0001-6556-8252]
\ead{moritz.diehl@imtek.uni-freiburg.de}
\credit{Conceptualization, Methodology, Formal analysis, Writing -- Review \& Editing, Supervision, Funding Acquisition}

\affiliation[1]{%
            organization={Department of Microsystems Engineering (IMTEK), University of Freiburg},
            addressline={Georges-Koehler-Allee 102}, 
            city={Freiburg},
            postcode={79110}, 
            country={Germany}}

\affiliation[2]{%
            organization={Department of Mathematics, University of Freiburg},
            addressline={Ernst-Zermelo-Str. 1}, 
            city={Freiburg},
            postcode={79114}, 
            country={Germany}}

\cortext[1]{Corresponding author}

\begin{abstract}
We present Riccati-ZORO, an algorithm for tube-based optimal control problems (OCP).
Tube OCPs predict a tube of trajectories to capture predictive uncertainty.
The tube induces a constraint tightening via additional backoff terms.
This backoff can significantly affect the performance, and thus implicitly defines a cost of uncertainty.
Optimizing the feedback law used to predict the tube can significantly reduce the backoffs, but its online computation is challenging.

Riccati-ZORO jointly optimizes the nominal trajectory and uncertainty tube based on a heuristic uncertainty cost design.
The algorithm alternates between two subproblems:
(i) a nominal OCP with fixed backoffs,
(ii) an unconstrained tube OCP, which optimizes the feedback gains for a fixed nominal trajectory.
For the tube optimization, we propose a cost function informed by the proximity of the nominal trajectory to constraints, prioritizing reduction of the corresponding backoffs.
These ideas are developed for ellipsoidal tubes under linear state feedback.
In this case, the decomposition into the two subproblems substantially reduces the computational complexity with respect to the state dimension from $\bigO(n_x^6)$ to $\bigO(n_x^3)$, i.e., the complexity of a nominal OCP.

We investigate the algorithm in numerical experiments, and provide two open-source implementations:
a prototyping version in CasADi and a high-performance implementation integrated into the \texttt{acados} OCP solver.
\end{abstract}




\maketitle

\section{Introduction}

Uncertainty-aware model predictive control (MPC), such as stochastic or robust MPC, aims for an explicit treatment of predictive uncertainty \citep{Rawlings2017,Kouvaritakis2016}.
Typically, the main motivation is to make the constraints robust against a given uncertainty model.

Existing approaches can be divided into two main families:
Scenario- or tree-based methods represent uncertainty discretely \citep{Scokaert1998,Calafiore2006,Lucia2014}, whereas tube-based methods use continuous uncertainty models \citep{Langson2004, Mayne2011,Rakovic2019, Koehler2021}. 
In tube-based methods, the uncertainty tube is typically predicted under the assumption of a simple feedback law, which mitigates the unrealistically large tubes resulting from an open-loop prediction.
In general, this feedback law is only used for the uncertainty prediction within the MPC problem. 
The actually applied feedback is determined by the MPC policy.

To ease the online computational burden, this internal feedback law is typically predesigned offline.
However, this comes at a loss of performance compared to online optimization, which enables the optimal shaping of the uncertainty tube for the currently planned trajectory \citep{Nagy2004, Messerer2021, Leeman2023}.
Furthermore, and especially for nonlinear systems, it is often unclear how to predesign a simple feedback law that performs adequately over all relevant state space regions.

\subsection{Problem structure}
In this paper, we consider tube-based optimal control problems (OCP) of the form \citep{Nagy2004, Gillis2015, Messerer2021, Leister2025, Belvedere2025, Zhang2025}
{\setlength{\mathindent}{7pt}%
\begin{mini!}[1]
   {x, u, K, P}
   {\displaystyle\sum_{k=0}^{N-1} l_k(x_k, u_k) +
   l_N(x_N)
   }
   { \label{ocp:rob} }{}
   \addConstraint{x_0}{=\bar x_0, \quad P_0 = \bar P_0, \label{ocp:rob_initial}}
   \addConstraint{x_{k+1}}{=f_{k}(x_k, u_k),\label{ocp:rob_dyn_nom}}{k=0,\dots,N-1,}
   \addConstraint{P_{k+1}}{ =  \psi_k(P_k, K_k, x_k, u_k), \label{ocp:rob_dyn_unc}}{k=0,\dots,N-1,}
   \addConstraint{0}{\geq h_k(x_k, u_k) \nonumber}
   \addConstraint{}{\quad + b_k(P_k, K_k, x_k, u_k),\label{ocp:rob:ineq_constr_stage}}{k=0,\dots,N-1,}
   \addConstraint{0}{\geq h_N(x_N) + b_N(P_N, x_N).\label{ocp:rob:ineq_constr_terminal}}
\end{mini!}
}

Here, $x = (x_0, \dots, x_{N})$ and $u = (u_0, \dots, u_{N-1})$ denote the nominal state and input trajectories, with nonlinear dynamics~$f_k$.
The uncertainty tube around the nominal trajectory is parametrized by $P = (P_0, \dots, P_{N})$, and propagated under a feedback policy parametrized by $K=(K_0, \dots, K_{N-1})$.
In this paper, we focus on the case where $P_k$ parametrizes an ellipsoidal tube, and $K_k$ parametrizes linear state feedback.

The uncertainty tube induces a backoff $b_k$, by which the nonlinear constraints $h_k$ on the nominal trajectory need to be tightened to stay safe.
This implicitly defines a cost of uncertainty through the impact of the constraint backoffs on the nominal trajectory.

Crucially, \eqref{ocp:rob} optimizes the feedback parameters $K_k$, such that it optimally handles uncertainty for the given tube and policy parametrization.
However, this is challenging to solve, both with respect to the computational cost and algorithmic robustness.

\subsection{Related work}

Problems of structure \eqref{ocp:rob} appear in robust and stochastic MPC for different tube and feedback parametrizations.
In this paper, we focus on the case of ellipsoidal tubes under linear state feedback, which are propagated based on a linearization at the nonlinear nominal trajectory \citep{Nagy2004, Diehl2006c, Houska2010, Messerer2021, Belvedere2025, Zhang2025}.
Linearization-based propagation incurs a linearization error, which may be conservatively overbounded \citep{Houska2011d, Koller2018, Kim2024b, Leister2025}.
The tube propagation can also be combined with learned uncertainty models such as Gaussian Processes (GP) \citep{Koller2018, Hewing2020}.
While state feedback assumes an exact measurement of the state at each time step, there are also extensions to the setting of output feedback \citep{Farina2015, Messerer2023}.

The optimization over state feedback in \eqref{ocp:rob} is highly nonlinear, which poses a challenge for general purpose solvers.
\citet{Goulart2006} show that linear feedback over all past states can be equivalently reformulated as linear disturbance feedback.
In contrast to state feedback, the optimization over disturbance feedback yields convex problems (for linear dynamics and convex costs and constraints).
However, this results in $\bigO(N^2)$ optimization variables, compared to $\bigO(N)$ for standard state feedback.
\citet{Leeman2025a} derive a disturbance feedback formulation for the nonlinear case, with linearization-based disturbance propagation and error bounds.

Several algorithms exploit the specific structure of \eqref{ocp:rob}.
For the case of ellipsoidal tubes, this reduces the complexity per iteration from $\bigO(n_x^6)$ to $\bigO(n_x^3)$ compared with a standard OCP solver.
\citet{Zanelli2021} propose ZORO (zero-order robust optimization) for the case of predesigned feedback gains.
As a zero-order method, ZORO converges in general to a feasible point in the neighborhood of an exact solution \citep{Zanelli2021}.
Convergence to the exact solution can be recovered by using an appropriate gradient correction \citep{Feng2020}.
A high-performance implementation of ZORO in \texttt{acados} \citep{Verschueren2021} is presented in \citet{Frey2024}. 
Real-time feasibility of ZORO for collision avoidance in mobile robotics has been demonstrated in real-world experiments \citep{Gao2023}.
\citet{Lahr2023} extend ZORO to GP-based MPC,
with an efficient implementation in \citet{Lahr2024} augmenting the \texttt{acados} ZORO implementation.
\citet{Zhang2024a} adapt ZORO for time-optimal motion planning.

SIRO (Sequential inexact robust optimization) \citep{Messerer2021} solves \eqref{ocp:rob} exactly in the special case of Lyapunov dynamics under linear state feedback.
Here, the optimality conditions of \eqref{ocp:rob} define the optimal feedback gains as the solution of a linear quadratic regulator (LQR) problem, in which the uncertainty in the constraint directions is weighted by the corresponding Lagrange multipliers.
Local convergence of SIRO is proved in \citet{Messerer2021} under the assumption of a stable active set.
However, subsequent experiments have revealed that for some practically relevant problems this assumption is not always fulfilled.
SIRO is extended to linear disturbance feedback in \citet{Leeman2024}, and adapted for time-optimal motion planning in \citet{Zhang2025}.

Furthermore, several algorithms perform heuristic online optimization of the feedback gains.
The algorithm in \citet{Leister2025} can be seen as a modification of ZORO, updating the feedback gains in each iteration by solving an LQR problem.
\citet{Belvedere2025} use the sensitivities of the nominal OCP solution as an approximation of future MPC feedback.
\citet{Kim2024b} solve a semidefinite program (SDP) to jointly update the tube and feedback gain parameters for the current nominal trajectory.

\subsection{Contribution and Outline}
In this paper, we propose the algorithm Riccati-ZORO for the heuristic optimization of feedback gains in \eqref{ocp:rob}.
This generalizes ZORO \citep{Zanelli2021, Frey2024}, which addresses \eqref{ocp:rob} for the case of predesigned feedback.
Like \citet{Leister2025}, Riccati-ZORO optimizes the feedback gains in each iteration with respect to a heuristically designed cost related to LQR.
In contrast to the constant cost weights in \citet{Leister2025}, Riccati-ZORO updates the weights based on constraint proximity to mimic the effect of exact feedback gain optimization.
From this perspective, Riccati-ZORO is a heuristic variation of SIRO \citep{Messerer2021}, giving the user more control over the tube shaping subproblem.

Riccati-ZORO is made available through two open-source implementations: (i) a prototyping version in CasADi \citep{Andersson2019}, (ii) a high-performance implementation in \texttt{acados} \citep{Verschueren2021}, extending the ZORO implementation \citep{Frey2024}.

This paper is structured as follows.
Section~\ref{sec:ellipsoidOCP} introduces the treated ellipsoidal tube OCP,
and Section~\ref{sec:tubeLQR} revisits the LQR from a tube perspective.
Section~\ref{sec:riccatizoro} defines the proposed algorithm and cost design, with a brief discussion of optimality and convergence in Section~\ref{sec:analysis}.
Section~\ref{sec:implementation} describes the implementation in \texttt{acados}, followed by numerical experiments in Section~\ref{sec:numerical} and a concluding Section~\ref{sec:conclusion}.

\subsection{Notation and preliminaries}
For vectors $x\in\R^n$, $y\in\R^m$ we denote by $(x,y) \defeq [x^\top, y^\top]^\top$ their vertical concatenation.
For a vector-valued function $f\colon\R^n \to \R^q$, we denote by $\nabla f(x) \in\R^{n \times q}$ the gradient, which is the transpose of the Jacobian $\dpartial{f(x)}{x}$. For a function with multiple arguments, $g\colon \R^n \times\R^m \to \R^q$, we use $\nabla g(x,y)$ to denote the vertically concatenated gradient w.r.t. both arguments, or specify the argument as $\nabla_x g(x,y)$.
Each symmetric positive semidefinite matrix $Q\in \SSS_+^n$ defines an ellipsoid $\calE(c, Q) \defeq  \{c + Q^{1/2} x \mid x \in \R^n, \,x^\top x \leq 1 \} $, with center $c \in \R^n$.
The Euclidean norm weighted by $Q$  is denoted as $\lVert x \rVert_Q = \sqrt{ x^\top Q x }$.
Using the cyclic property of the trace operator $\tr(\cdot)$, the trace trick allows the reformulation
\begin{equation}\label{eq:tracetrick}
   x^\top Q x = \tr( x^\top Q x) = \tr(Q x x^\top).
\end{equation}

\section{Linearization-based ellipsoidal tube OCP}
\label{sec:ellipsoidOCP}
We now derive a problem of structure \eqref{ocp:rob} for the case of ellipsoidal uncertainty sets, which are propagated based on a linearization along the nominal trajectory \citep{Messerer2021, Diehl2006c, Houska2010}.
A linear state feedback law is used to counteract tube growth.

Consider the uncertain nonlinear dynamics
\begin{equation} \label{eq:dynamical_system}
   \tilde x_0 = \bar x_0, \;\;\; \tilde x_{k+1} = \tilde f_k(\tilde x_k, \tilde u_k, w_k),\;\;\; k=0,\dots, N-1,
\end{equation}
with uncertain state $\tilde x_k\in\R^{n_x}$, control $\tilde u_k\in\R^{n_u}$ and perturbation $w_k \in \R^{n_w}$. 
We assume that the perturbation trajectory jointly lies inside an ellipsoidal set, $w = (w_0, \dots, w_{N-1}) \in \calE(0, W)$, with $W=\blkdiag(W_0, \dots, W_{N-1})$.

We describe the resulting uncertainty of the state trajectory by ellipsoidal tubes $\calE(x_k, P_k)$, with $P_k \in \SSS^{n_x}_+$.
Propagating the system for $w_k = 0$ yields the nominal trajectory,
\begin{equation}
   x_0 = \bar x_0, \;\;\;  x_{k+1} = \tilde f_k(x_k, u_k, 0),\;\;\; k=0,\dots, N-1.
\end{equation}
The linear state feedback law
\begin{equation} \label{eq:linearstatefeedback}
   \tilde u_k =  u_k + K_k (\tilde x_k - x_k),\quad\text{with}\;\; K_k \in \R^{n_u \times n_x},
\end{equation}
is used to counteract tube growth by reacting to deviations from the nominal trajectory.
Based on a linearization at the nominal trajectory, we obtain the ellipsoid dynamics as
\begin{equation} \label{eq:ellipsoid_dynamics}
   \begin{aligned}
      P_{k+1}&= (A_k + B_k K_k) P_k (A_k + B_k K_k)^\top + \overbrace{\Gamma_k W_k \Gamma_k^\top}^{\eqdef \tilde W_k} \\ 
       &\eqdef  \psi_k(x_k, u_k, P_k, K_k),
   \end{aligned}
\end{equation}
where the sensitivity matrices $A_k$, $B_k$, $\Gamma_k$ depend on the nominal trajectory as
\begin{equation}
\begin{aligned}
      A_k &= \nabla_x f_k(x_k, u_k, 0)^\top, \quad 
      B_k = \nabla_u f_k(x_k, u_k, 0)^\top, \\
      \Gamma_k &= \nabla_w f_k(x_k, u_k, 0)^\top, \quad k=0,\dots, N-1.
\end{aligned}
\end{equation}

Given the ellipsoidal tube $\tilde x_k \in \calE(x_k, P_k)$, linear state feedback \eqref{eq:linearstatefeedback}, and a linearization of the constraint functions $h_k^i(\tilde x_k, \tilde u_k)$ at the nominal trajectory, the backoffs in the robustified constraints \eqref{ocp:rob:ineq_constr_stage}, \eqref{ocp:rob:ineq_constr_terminal} can be obtained as the worst-case deviation from the nominal constraint value,
\begin{subequations} \label{eq:backoff}
   \begin{flalign}
      b_k^i(x_k, u_k, P_k, K_k) &= \big\lVert \nabla_x h_k^i(x_k, u_k) + K_k^\top \nabla_u h_k^i(x_k, u_k)  \big\rVert_{P_k},&
      \\
         b_N^i(x_N,P_N) &= \big\lVert \nabla_x h_N^i(x_N) \big\rVert_{P_N}, &
\end{flalign}
\end{subequations}
for $i=1,\dots, n_h$, $k=0,\dots,N-1$.

\begin{remark}[Linearization with respect to uncertainty]
Under the given assumptions on $w$, the ellipsoid dynamics~\eqref{eq:ellipsoid_dynamics} and constraint robustification \eqref{eq:backoff} are exact for linear dynamics $f_k$ and constraints $h_k$, but otherwise subject to a linearization error.
This error can be accounted for by adding additional terms in \eqref{eq:ellipsoid_dynamics} and \eqref{eq:backoff}, cf., e.g., \citet{Houska2011d, Kim2024b, Leister2025, Leeman2025a}.
\end{remark}

\begin{remark}[Complexity of the standard OCP structure] \label{rem:complexity}
Tube OCP \eqref{ocp:rob} has a standard OCP sparsity structure with respect to the augmented state $\breve x_k = (x_k, \vecop(P_k))$ and augmented control $\breve u_k = (u_k, \vecop(K_k))$.
If $P_k$, $K_k$ parametrize an ellipsoidal tube under linear state feedback, the augmented state and control have dimensions $\breve n_x = n_x + \tfrac{1}{2}(n_x + n_x^2)$ and $\breve n_u = n_u + n_x n_u$, where for $\vecop(P_k)$ the symmetry of $P_k$ is exploited.
When solving \eqref{ocp:rob} with a standard OCP algorithm, the computational complexity is
$\bigO((\breve n_x^3 + \breve n_x^2\breve n_u + \breve n_x \breve n_u^2 + \breve n_u^3)N)$ per iteration \citep{Rawlings2017}.
With respect to the original state and control dimension, this corresponds to
$\bigO((n_x^6 + n_x^5 n_u + n_x^4 n_u^2 + n_x^3 n_u^3)N)$.
\end{remark}

\begin{remark}[Stochastic interpretation]\label{rem:stoch}
   In this section, we have derived an OCP of structure \eqref{ocp:rob} under the assumption of a set-based bound on the perturbation trajectory.
   The same OCP can be derived in a stochastic setting, for independent noise $w_k$ with zero mean and covariance $W_k$.
   In this case, $x_k$, $P_k$ parametrize a distribution in state space with mean $x_k$ and covariance $P_k$.
   The constraints \eqref{ocp:rob:ineq_constr_stage}, \eqref{ocp:rob:ineq_constr_terminal} are chance constraints, with the backoffs $b_k^i$ corresponding to the standard deviation in constraint direction.
   To achieve a desired level of constraint satisfaction probability, the backoff needs to be multiplied by a corresponding factor \citep{Mesbah2016}.
\end{remark}

\section{Linear quadratic regulator in tube space}
\label{sec:tubeLQR}

A central building block of the proposed algorithm is the update of the feedback gains and ellipsoidal tube for a given nominal trajectory.

Recall that $P_k$ describes the deviations from the nominal trajectory as an ellipsoidal set, $\tilde x_k - x_k \in \calE(0, P_k)$.
Linear feedback on these deviations \eqref{eq:linearstatefeedback} induces a joint ellipsoid in the combined state-control space:
\begin{equation}
   \bmat \tilde x_k - x_k \\ \tilde u_k - u_k \emat \in \calE(0, \bmat P_k&P_k K_k^\top \\ K_k P_k & K_k P_k K_k^\top  \emat).
\end{equation}

We can optimize the trajectory of this ellipsoid, balancing uncertainty in different directions, by minimizing its weighted trace,

\begin{samepage}
\begin{subequations}\label{ocp:lqr_matrixform}
    \begin{flalign}
        \min_{\displaystyle K, P} \qquad &
        \sum_{k=0}^{N-1} \tr\Bigg(\!\overbrace{\bmat Q_k & S_k^\top \\ S_k & R_k \emat}^{C_k} \overbrace{\bmat P_k&P_k K_k^\top \\ K_k P_k & K_k P_k K_k^\top \emat}^{\check P_k} \! \Bigg)  \nonumber &&\\
        &  \;\;\; + \tr (Q_N P_N) \label{ocp:lqr_matrixform:objective}
        \\
        \mathrm{s.t.}\quad \;
        P_0&=\bar P_0,\\
        P_{k+1} &= (A_k + B_k K_k) P_k (A_k + B_k K_k)^\top + \tilde W_k, \nonumber &&\\
        &\qquad\qquad\qquad\qquad k=0,\dots, N-1. \label{ocp:lqr_matrixform:dynamics}
    \end{flalign}
\end{subequations}
\end{samepage}

This is an OCP in ellipsoidal tube space with state $P_k$ and control $K_k$.
To interpret the objective, we note that the trace of a positive semidefinite matrix is given by the sum of its eigenvalues, which in turn are the squared axis lengths of the corresponding ellipsoid.
In this sense, \eqref{ocp:lqr_matrixform} minimizes squared deviations in the state-control space, weighted by $C_k$, such that we also refer to $C_k$ as a ``Hessian''.

In fact, \eqref{ocp:lqr_matrixform} can alternatively be derived as a reformulation of the standard stochastic LQR problem, cf. \citet{Messerer2021}.
In this context, $P_k$ and $\check P_k$ are covariance matrices, and \eqref{ocp:lqr_matrixform} minimizes the expected value of the squared deviation.
This follows from
\begin{equation}\label{eq:exv_of_quadratic}
   \begin{aligned}
      &\EXV_{\tilde z_k} \{ (\tilde z_k - z_k)^\top C_k  (\tilde z_k - z_k) \}\\
      &=\EXV_{\tilde z_k} \{ \tr\left(C_k  (\tilde z_k - z_k) (\tilde z_k - z_k)^\top\right) \} = \tr (C_k \check P_k),
   \end{aligned}
\end{equation}
where we used $\tilde z_k = (\tilde x_k, \tilde u_k)$ for the uncertain state-control trajectory with mean $ z_k$ and covariance $\check P_k$, and the second line follows from the trace trick \eqref{eq:tracetrick}.

As a classical result, for the stochastic LQR the optimal feedback gains are given by a Riccati recursion, which, due to certainty equivalence, is independent of the noise contribution~$\tilde W$ \citep{Stengel1986}.
This holds true for the algebraic structure of \eqref{ocp:lqr_matrixform} in general, independent of the interpretation of $P_k$ as covariances or ellipsoidal sets.

\begin{lemma}[\citet{Messerer2021}]\label{lem:riccati_for_matrixlqr}
   Let $Q_N \succeq 0$, $C_k \succeq 0$, $Q_k \succeq 0$, $R_k \succ 0$, for $k=0, \dots, N-1$. Further, let $P_k \succ 0$ for all $P_k$ reachable with dynamics \eqref{ocp:lqr_matrixform:dynamics}, $k=0,\dots,N$. Then, the solution to \eqref{ocp:lqr_matrixform} is uniquely defined by the backward Riccati recursion
   \begin{subequations}\label{eq:riccati}
      \begin{flalign}
      V_N &= Q_N,&\\
      K_k^\star &= - (R_k + B_k^\top V_{k+1} B_k)^{\shortminus1} (S_k + B_k^\top V_{k+1} A_k),&\\
      V_k &= Q_k + A_k^\top V_{k+1} A_k + (S_k^\top + A_k^\top V_{k+1} B_k)K_k^\star,&
   \end{flalign}
\end{subequations}
   for $k = N-1, \dots, 0$, followed, for $k=0, \dots, N-1$, by the forward Lyapunov recursion
\begin{subequations} \label{eq:lyapunovForward}
         \begin{flalign} 
         P_0^\star &= \bar P_0, \\
         P_{k+1}^\star &= (A_k + B_k K_k^\star) P_k^\star (A_k + B_k K_k^\star)^\top + \tilde W_k.&
      \end{flalign}
\end{subequations}
\end{lemma}

\section{Riccati-ZORO}
\label{sec:riccatizoro}

\begin{figure}
   \includegraphics[width=\columnwidth]{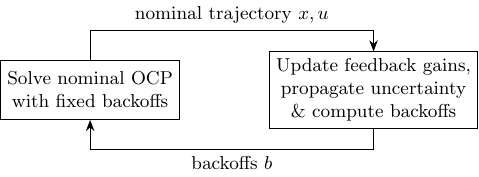}
   \centering
   \caption{Schematic structure of Riccati-ZORO}
   \label{fig:algo_schematic}
\end{figure}

In this section we describe the algorithm Riccati-ZORO, which targets problems of structure \eqref{ocp:rob}, but optimizes the feedback gains based on a heuristic cost instead of exactly.
We consider two choices for the heuristic cost, and discuss the relation to the algorithms ZORO and SIRO. 

\subsection{Structure of the algorithm}

Riccati-ZORO alternates between solving a nominal OCP with fixed backoffs $b$, and a feedback gain optimization problem for a fixed nominal trajectory, which subsequently yields updated backoffs.
This structure is visualized in Fig.~\ref{fig:algo_schematic}.

\subsubsection{Nominal OCP with fixed backoffs}
By computing the backoffs for fixed values of $P$ and $K$ in \eqref{ocp:rob}, we obtain the first subproblem as
{\setlength{\mathindent}{20pt}%
\begin{mini!}[1]
   {x, u}
   {\sum_{k=0}^{N-1} l_k(x_k, u_k) + l_N(x_N)}
   {\label{ocp:perturbed_nominal}}{}
   \addConstraint{x_0}{=\bar x_0,}
   \addConstraint{x_{k+1}}{=f_k(x_k, u_k), \label{ocp:perturbed_nominal:ineq_stage} }{k=0, \dots, N-1,}
   \addConstraint{0}{\geq h_k(x_k, u_k) + b_k,\quad }{k=0, \dots, N-1,}
   \addConstraint{0}{\geq h_N(x_N) + b_N,}
\end{mini!}
}

which has the structure of a standard nominal OCP.

The solution of \eqref{ocp:perturbed_nominal} yields an updated nominal trajectory.
The second subproblem, described below, optimizes the feedback gains for this updated trajectory by solving \eqref{ocp:lqr_matrixform}.
From the resulting tube, updated backoffs are computed, cf. Fig.~\ref{fig:algo_schematic}.

Treating the backoffs as fixed neglects their dependence on the tube variables, which depend on the nominal trajectory.
Since this disregards first-order derivatives, the resulting algorithm is a zero-order method, as discussed in Section~\ref{sec:analysis}.

For intermediate iterations of Riccati-ZORO, subproblem \eqref{ocp:perturbed_nominal} can be infeasible due to large backoffs, even if \eqref{ocp:rob} is feasible.
In this case, one can enforce the inequality constraint using an exact penalty on constraint violation, instead of explicit hard constraints \citep{Nocedal2006}.
This yields an exact solution of \eqref{ocp:perturbed_nominal} if feasible, and otherwise minimizes constraint violation.

For simplicity of presentation, we assume that in each iteration, \eqref{ocp:perturbed_nominal} is solved up to convergence.
Alternatively, one can update the backoffs sooner, e.g., after performing a single SQP iteration on \eqref{ocp:perturbed_nominal}, thereby avoiding unnecessary iterations on a problem that is ultimately not of interest. 
This is discussed in more detail in Section~\ref{sec:implementation}.

\subsubsection{Tube optimization by Riccati recursions}
In the second subproblem, the feedback gains are updated for the current nominal trajectory.
This is achieved by solving \eqref{ocp:lqr_matrixform}, which optimizes the trajectory of the joint state-control ellipsoid by minimizing its size, using the weighted trace as a metric.
Crucially, we allow the weighting matrix $C_k$ to depend on the nominal trajectory.
By Lemma~\ref{lem:riccati_for_matrixlqr}, the solution of \eqref{ocp:lqr_matrixform} is given by a Riccati recursion followed by the ellipsoid propagation.

Overall, the computational steps of Riccati-ZORO are given in Algorithm~\ref{algo:Riccati-ZORO}.

\begin{algorithm}[t]
   \caption{Riccati-ZORO}
   \label{algo:Riccati-ZORO}
   \begin{algorithmic}
      \State \textbf{Input:} Initial guess $x$, $u$
      \textit{(e.g. solution to nominal OCP)}
      \While{\texttt{not\_converged}}
         \State $K \gets$ \texttt{riccati\_recursion($x, u$)}
            \Comment{\eqref{eq:riccati}}
         \State $P \gets$ \texttt{propagate\_ellipsoids($x, u, K$)}
            \Comment{\eqref{eq:lyapunovForward}}
         \State $b \gets$ \texttt{compute\_backoffs($x, u, P, K$)}
            \Comment{\eqref{eq:backoff}}
         \State $x, u \gets $ \texttt{solve\_perturbed\_ocp($b$)}
            \Comment{\eqref{ocp:perturbed_nominal}}
      \EndWhile
      \State \textbf{return:} $x, u, P, K$
 \end{algorithmic}
\end{algorithm}

\subsection{Uncertainty cost design}
We consider two choices for the uncertainty weighting matrices $C_k$ used in the tube optimization subproblem \eqref{ocp:lqr_matrixform}.

\subsubsection{Trajectory-independent weighting}
The most straightforward choice is to update the feedback gains by solving \eqref{ocp:lqr_matrixform} with cost matrices independent of the nominal trajectory, as proposed in \citet{Leister2025}.
This provides a trajectory-independent balancing of uncertainty across all directions of the state-control space, but does not systematically trade off uncertainties between directions according to their performance impact.
For tuning this cost, the same intuitions as for standard LQR hold.

\subsubsection{Constraint-adaptive weighting}
Here, we derive a weighting scheme that prioritizes reducing the uncertainty in constraint directions, yielding smaller backoffs, based on the proximity of the nominal trajectory to the constraint.
This weighting is most intuitively derived through the approximate expected value of a logarithmic barrier function, though as a heuristic, it remains independent of the specific interpretation of the ellipsoids.

First, consider a nonlinear cost function $\hat l_k(\tilde z_k)$.
Given uncertain $\tilde z_k$ with mean $z_k$ and covariance $\check P_k$, we can approximate its expectation using a second-order Taylor expansion as
\begin{equation}\label{eq:expectation_over_quadratic_approx}
      \begin{aligned}
         \EXV_{\tilde z_k} \{ \hat l_k(z_k)\} 
         &\approx \EXV_{\tilde z_k} \{ \hat l_k(z_k) + \nabla   \hat l_k(z_k)^\top(\tilde z_k-z_k) \\
         & \quad + \tfrac{1}{2} (\tilde z_k - z_k)^\top \nabla^2   \hat l_k(z_k) (\tilde z_k -  z_k) \} \\
         &= \hat l_k(z_k) +  \tfrac{1}{2}\tr(\nabla^2  \hat l_k(z_k) \check P_k),
      \end{aligned}
\end{equation}
where the second-order term follows from \eqref{eq:exv_of_quadratic}.
The first term is the nominal cost, whereas the second term penalizes covariance weighted by the nominal curvature.

Now, consider enforcing the inequality $h_k^i(\tilde z_k) \leq 0$ with $\tilde z_k = (\tilde x_k, \tilde u_k)$ with a log-barrier term $\hat l_k^i(\tilde z_k) = \tau_k^i \phi(h_k^i(\tilde z_k))$,
with logarithmic barrier function $\phi(\eta) = - \log(-\eta)$ weighted by $\tau_k^i > 0$.
We approximate its expectation based on \eqref{eq:expectation_over_quadratic_approx}, using the approximate Hessian
\begin{equation}
   \nabla^2  \hat l_k^i(z_k) 
   \approx
   \tau_k^i
   \nabla  h_k^i(z_k) \phi''(h_k^i(z_k))\nabla h_k^i(z_k)^\top \succeq 0,
\end{equation}
which results from an additional linearization of the constraint function.
Since here we are only interested in deriving a heuristic cost on uncertainty, we drop the nominal term in \eqref{eq:expectation_over_quadratic_approx}.
Adding an additional regularization term $\bar C_k \succeq 0$ yields an objective of structure \eqref{ocp:lqr_matrixform:objective} with cost matrices
\begin{equation}\label{eq:costmat_hessbarrier}
   C_k(z_k) = \bar C_k + \sum_{i=0}^{n_{h_k}} \tau_k^i  \phi''(h_k^i(z_k))  \nabla h_k^i(z_k) \nabla h_k^i(z_k)^\top.
\end{equation}

\begin{figure}
   \centering
   \includegraphics[width=\columnwidth]{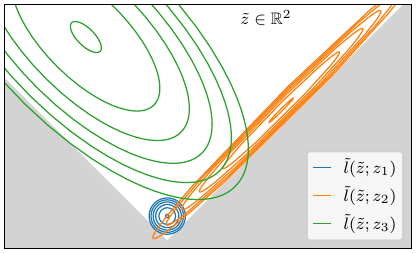}
   \caption{Visualization of the constraint-adaptive uncertainty weighting $C(z)$ in \eqref{eq:costmat_hessbarrier}, for $z \in \R^2$, and with $\bar C = 0$ and $\tau = 1$.
   The figure shows the contour lines of the quadratic uncertainty penalization $\tilde l(\tilde z; z) = (\tilde z - z)^\top C(z) (\tilde z - z)$ for the three exemplary values $z_1, z_2, z_3$.
   The gray area is the infeasible region.
   As $z$ approaches the boundaries of the feasible set, the penalty becomes steeper for the corresponding constraint direction.%
   }
   \label{fig:barrier_hess}
\end{figure}

The terms inside the sum correspond to the squared backoff terms penalized by the barrier curvature, as can be seen from the backoff definition \eqref{eq:backoff} and the trace trick \eqref{eq:tracetrick}.
This smoothly increases the penalty on a backoff term as the nominal trajectory approaches the corresponding constraint boundary.
The corresponding quadratic cost on deviations from the nominal trajectory, cf.~\eqref{eq:exv_of_quadratic}, is visualized in Fig.~\ref{fig:barrier_hess}.

By adapting the weights $\tau_k^i$, uncertainty can be traded off across different constraints.
A small trajectory-independent term $\bar C_k$ can serve as a baseline. 

As a heuristic, the performance of Riccati-ZORO depends on the uncertainty cost design.
If a good approximation of the optimal solution of \eqref{ocp:rob} can be achieved by keeping all backoffs small, Riccati-ZORO is expected to perform well with minor tuning effort.
However, if performance is highly sensitive to accurately capturing the optimal balance between backoffs, closely approximating the optimal solution of \eqref{ocp:rob} using Riccati-ZORO may require significant tuning effort.

\subsection{Computational complexity}

Solving the nominal OCP \eqref{ocp:perturbed_nominal} with standard OCP algorithms incurs a computational cost of 
$\bigO((n_x^3 + n_x^2n_u + n_x n_u^2 + n_u^3)N)$ \citep{Rawlings2017}.
The Riccati recursion \eqref{eq:riccati} has the same complexity, and the Lyapunov recursion \eqref{eq:lyapunovForward} has $\bigO((n_x^3 + n_x^2n_u)N)$.
Thus, the overall complexity of Riccati-ZORO is $\bigO((n_x^3 + n_x^2n_u + n_x n_u^2 + n_u^3)N)$ per iteration, which is the same complexity as a nominal OCP.
In contrast, solving \eqref{ocp:rob} directly with standard OCP solvers incurs a cost of $\bigO((n_x^6 + n_x^5 n_u + n_x^4 n_u^2 + n_x^3 n_u^3)N)$, cf. Remark~\ref{rem:complexity}.

\subsection{Relation to ZORO and SIRO}
Both ZORO and SIRO are closely related to Riccati-ZORO.

ZORO targets a simplified version of \eqref{ocp:rob}, in which the feedback gains $K$ are set to a predetermined value~$\bar K$ \citep{Zanelli2021, Frey2024}.
Thus, in each iteration it only updates the ellipsoid propagation and the backoffs.
In this sense, Riccati-ZORO generalizes ZORO.
In both variants, the backoff values are essentially functions of the nominal trajectory only, since the tube parameters could in principle be eliminated. 
Thus, at a high level, they address an identical problem structure, such that the theoretical properties of ZORO \citep{Zanelli2021} can be generalized to Riccati-ZORO, as discussed in Section~\ref{sec:analysis}.

SIRO targets the full problem \eqref{ocp:rob} with uncertainty dynamics and backoffs in exactly the form \eqref{eq:ellipsoid_dynamics} and \eqref{eq:backoff}.
Like Riccati-ZORO, it updates the feedback parameters by solving a problem of structure \eqref{ocp:lqr_matrixform} in each iteration.
Crucially, in SIRO the cost matrices for this update are not chosen heuristically, but are derived from the optimality conditions of \eqref{ocp:rob}.
This results in
\begin{align}\label{eq:costmat_siro}
   C_k(z_k, \mu_k) &= \bar C_k + \sum_{i=0}^{n_{h_k}} \frac{\mu_k^i}{2 b_k^i}  \nabla h_k^i(z_k) \nabla h_k^i(z_k)^\top,
\end{align}
with $\mu_k^i$ the Lagrange multiplier of the constraint $h_k^i$ in \eqref{ocp:perturbed_nominal}.
These cost matrices exactly specify the true cost of each backoff.
However, the dependence of \eqref{eq:costmat_siro} on the multipliers $\mu_k^i$ can create an algorithmic instability, keeping SIRO from converging.
This relates to active set changes of \eqref{ocp:perturbed_nominal}, which cause a nonsmooth change in the value of $\mu_k^i$.
Given the structural similarity, the barrier Hessian in \eqref{eq:costmat_hessbarrier} can be understood as a heuristic smoothing of the exact cost \eqref{eq:costmat_siro}.
As a minor difference, SIRO as presented in \citet{Messerer2021} uses a gradient correction and is therefore not a zero-order method.
An overview of the three algorithms is given in Table~\ref{tab:comp_algs}.

\newcommand{\ndashes}[1]{\ifnum#1>0\kern0.8pt$-$\kern0.8pt\ndashes{\numexpr#1-1\relax}\fi}
\begin{table*}
   \caption{Comparison of the three related algorithms ZORO, Riccati-ZORO and SIRO.}
   \label{tab:comp_algs}
   \centering
   \begin{tabular}{@{}llll@{}}
      \toprule
                           & ZORO \citep{Zanelli2021}   &  Riccati-ZORO & SIRO \citep{Messerer2021} \\
      \midrule
   Online internal feedback optim.  & \xmark & \cmark{} (heuristic) & \cmark{} (optimal)   \\
   Backoff update  &  Lyapunov & Riccati \& Lyapunov &   Riccati \&  Lyapunov \\   
   Backoff update dependency  & nominal trajectory & nominal trajectory &  nominal traj. \& constr. multipliers \\   
   Active set changes & no impact on backoff &  smoothed by heuristic & can cause convergence issues \\
   Stationary points & feasible, suboptimal & feasible, suboptimal &  optimal \\
   Implementation in \texttt{acados} & \cmark{} \citep{Frey2024} & \cmark{} (this work) &  \xmark \\
      Comp. complexity per iteration
   & \multicolumn{3}{l@{}}{\kern-.85pt$|$\kern-2.8pt\ndashes{16}\kern2.2pt $\bigO((n_x^3 + n_x^2n_u + n_x n_u^2 + n_u^3)N)$\kern2.2pt\ndashes{16}\kern-2.8pt$|$\kern-.85pt} \\
   \bottomrule
   \end{tabular}
\end{table*}

\section{Analysis of optimality and convergence}
\label{sec:analysis}

In this section, we analyze the points to which Riccati-ZORO converges.
In comparison to the exact feedback gain optimization problem \eqref{ocp:rob}, Riccati-ZORO introduces two sources of suboptimality.
The first results from the heuristic choice of feedback gains; the second from the zero-order backoff updates.
In both cases, constraint satisfaction and any associated guarantees remain uncompromised.

We start by writing~\eqref{ocp:rob} at a more abstract level.
By eliminating the nominal and uncertainty state trajectories $x$ and $P$ based on their dynamics \eqref{ocp:rob_dyn_nom}, \eqref{ocp:rob_dyn_unc}, and correspondingly redefining the objective and inequality constraints, we obtain
\begin{equation} \label{nlp:high}
   \min_{\displaystyle u, K} \; f(u) \quad\mathrm{s.t.}\quad h(u) + \sigma  b(\bar u, K) \leq 0
\end{equation}
as the sequential formulation of \eqref{ocp:rob}.
We additionally introduce the parameter $\sigma > 0$ as the level of uncertainty. 
This corresponds to substituting $\sigma^2 \bar P_0$ in \eqref{ocp:rob_initial} and $\sigma^2 \tilde W_k$ in \eqref{eq:ellipsoid_dynamics}.

Instead of solving \eqref{nlp:high} exactly, Riccati-ZORO heuristically chooses the feedback gains~$K$ as the solution of an optimization problem parametrized by $u$,
\begin{equation}
  \tilde K(u) = \argmin_{\displaystyle K} \; F(u, K),
\end{equation}
which in our case is given by \eqref{ocp:lqr_matrixform}.
This heuristic choice of feedback gains corresponds to a restriction of \eqref{nlp:high} to
\begin{equation} \label{nlp:high_m_of_y}
   \min_{\displaystyle u} \; f(u) \quad\mathrm{s.t.}\quad h(u) + \sigma  b( u, \tilde K( u)) \leq 0,
\end{equation}
where the nominal trajectory remains the only degree of freedom.
In general, a heuristic choice of $\tilde K(u)$ will not exactly account for the implicit cost of the backoffs, such that solving \eqref{nlp:high_m_of_y} yields a suboptimal solution of \eqref{nlp:high}.

Iteratively solving \eqref{nlp:high_m_of_y} for fixed backoffs neglects their gradient.
This corresponds to a zero-order approximation of the constraint gradient \citep{Zanelli2021}.
Consequently, stationary points $\hat u$ of Riccati-ZORO are not exact solutions of \eqref{nlp:high_m_of_y}.
Instead, they are solutions of a perturbed version of \eqref{nlp:high_m_of_y}, given by
\begin{equation} \label{nlp:high_m_of_y_perturbed}
   \min_{\displaystyle u} \; f(u) + \sigma \hat c^\top u  \quad\mathrm{s.t.}\quad h(u) + \sigma  b( u, \tilde K( u)) \leq 0,
\end{equation}
with gradient perturbation $\hat c = \dtotaltfrac{}{\hat u} b(\hat u,\tilde K(\hat u))^\top \hat \mu$ resulting from the neglected backoff gradient evaluated at the stationary point $\hat u$ with Lagrange multiplier $\hat \mu$.
Note that this perturbation only modifies the objective function.
The constraint satisfaction remains unaffected.

For more detail, we refer to \citet{Zanelli2021}, which analyzes the corresponding suboptimality for ZORO and proves local convergence of ZORO for sufficiently small~$\sigma$.
ZORO addresses \eqref{nlp:high_m_of_y} for the special case of a constant feedback gain, $\tilde K(u) \equiv \bar K$.
Thus, the only difference is that Riccati-ZORO additionally neglects the gradient of $\tilde K(u)$, such that the core arguments of the proofs in \citet{Zanelli2021} also apply to Riccati-ZORO.
\citet{Feng2020} recover optimality by applying a gradient correction, which in our case would need derivatives of the Riccati recursion and ellipsoid propagation.

\section{Efficient implementation in \texttt{\textup{acados}}}
\label{sec:implementation}

We provide an efficient open-source implementation of Riccati-ZORO in \texttt{acados} \citep{Verschueren2021} by extending the ZORO implementation described in \citet{Frey2024}.
The implementation is written in \texttt{C}, using \texttt{BLASFEO} \citep{Frison2018} for the linear algebra, and is accessible through the \texttt{acados} Python and Matlab interfaces.
The robust OCP can be conveniently defined on top of a standard \texttt{acados} OCP by specifying the uncertainty-related quantities.
Additionally, the online adaptation of externally computed covariance matrices, as needed, for example, in GP-MPC~\citep{Lahr2023}, is supported through the extension in \citet{Lahr2024}.

The solution of a nominal OCP with an SQP method requires at each time step the computation of the sensitivities of the dynamics and constraints, which is typically a major part of the computational load.
These are the same sensitivities used in the Riccati recursion \eqref{eq:riccati}, ellipsoid propagation \eqref{eq:lyapunovForward}, and backoff computation \eqref{eq:backoff}.

The \texttt{acados} custom update functionality introduced by \citet{Frey2024} allows executing a custom \texttt{C} function between the OCP solver calls, which can access and modify the solver data.
Riccati-ZORO extends the ZORO custom update. Overall, with each call it:
(i)~accesses the dynamic and constraint sensitivities, which have already been computed as part of the OCP solver call, (ii)~uses \texttt{BLASFEO} to perform the backoff update including Riccati recursion, ellipsoid propagation, and backoff computation, and
(iii)~updates the inequality constraint bounds with the new backoffs.

The implementation is released in \url{https://github.com/acados/acados/releases/tag/v0.5.2}.

\section{Numerical Experiments} \label{sec:numerical}

In this section, we compare Riccati-ZORO to exact feedback optimization, based on an implementation in CasADi \citep{Andersson2019} with IPOPT \citep{Waechter2006} as solver, followed by a benchmark of the computational performance of the \texttt{acados} implementation, and a case study from industrial research.
All experiments are performed on a standard laptop with an Intel i7-8565U CPU and 16 GB of RAM.
The code is available at \url{https://github.com/fmesserer/riccati-zoro}.

\begin{figure*}[t]
   \centering
   \includegraphics[width=\textwidth]{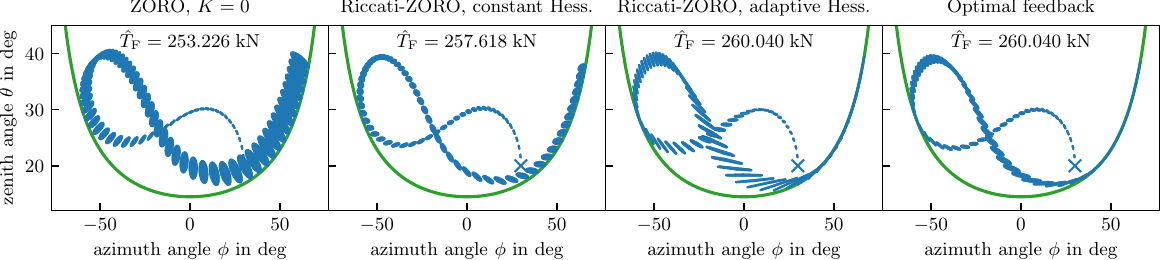}
   \caption{Stationary points of different algorithms for the towing kite problem, with the objective of maximizing the average thrusting force $\hat T_\mathrm{F}$.
   In comparison, the exact solution of the fixed feedback problem with $K=0$ with $\hat T_{\mathrm{F}} = 253.235~\mathrm{kN}$ looks indistinguishable from the ZORO solution. The nominal problem with zero backoff has $\hat T_{\mathrm{F}} = 260.086~\mathrm{kN}$.
   }
   \label{fig:kite}
\end{figure*}

\subsection{Towing Kite (with CasADi \& IPOPT)}
We use the towing kite example from \citep{Messerer2021} to compare the two Riccati-ZORO variants to the optimal solution of \eqref{ocp:rob}.
The towing kite is connected to a ship by a tether of constant length, and the objective is to maximize the average thrusting force, while respecting a minimum height constraint. 
The kite position is described by the angles $\theta$ and $\phi$ and its orientation by the angle $\psi$, resulting in the state $x=(\phi, \theta, \psi)$, with the steering deflection $u\in\R$ as control input.
The main source of uncertainty is the unknown wind speed. For more details, see \citep{Messerer2021}.
Due to its nonlinear economic objective function and highly nonlinear dynamics, this is a challenging problem.

We compare the constant-weight and constraint-adaptive variants of Riccati-ZORO to the exact solution obtained by SIRO and the solution obtained by ZORO with fixed feedback $K=0$.
For the constant Hessian variant of Riccati-ZORO, we use $Q_k=\eye$, $R_k = 10^{-2}$.
The results are shown in Fig.~\ref{fig:kite}.

With the feedback gains computed based on a constant Hessian, Riccati-ZORO achieves a balanced  uncertainty reduction in all directions.
In contrast, using the adaptive Hessian reduces uncertainty specifically in the constraint direction, such that for this example there is effectively no suboptimality with respect to exact feedback gain optimization.
For a rough comparison of timings, we initialize all solvers at the nominal trajectory.
ZORO, both Riccati-ZORO variants, and SIRO all need 6 or 7 iterations and $0.45\pm0.05\;\mathrm{s}$ to converge. 
In comparison, when solving \eqref{ocp:rob} directly, IPOPT needs $16.3\;\mathrm{s}$ for the problem with fixed feedback gains $K=0$, and does not converge for the full feedback gain optimization problem.

\subsection{Hanging chain of masses (with \texttt{acados})}
\begin{figure}[t]
   \centering
   \includegraphics[width=\columnwidth]{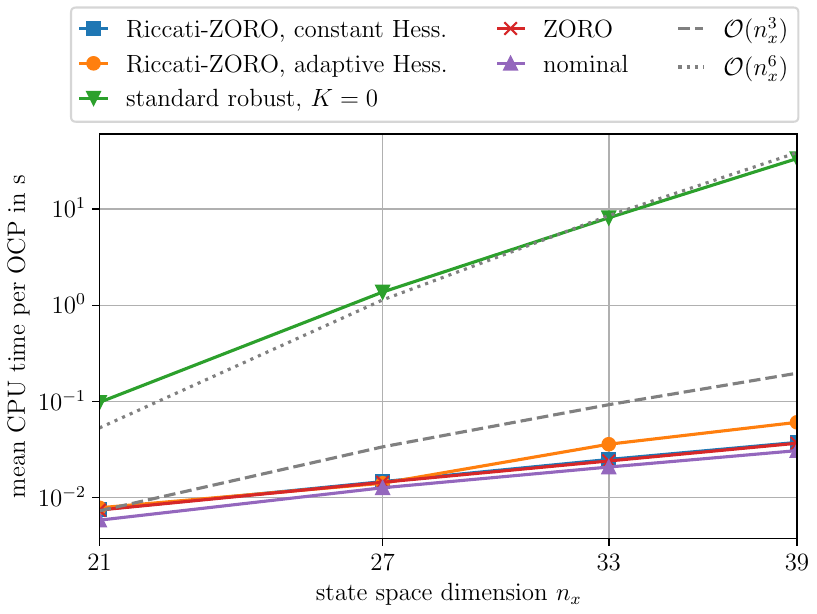}
   \caption{Computation times for the hanging chain of masses benchmark.}
   \label{fig:timings_order_nx}
\end{figure}

To investigate the complexity with respect to the state dimension, we use the hanging chain of masses benchmark in the setting detailed in \citet{Frey2024}.
By changing the number of masses in the chain, this example allows flexible scaling of the state space dimension $n_x$.
We compare the \texttt{acados} implementation of the two variants of Riccati-ZORO to nominal MPC and ZORO with fixed feedback $K=0$.
As a reference, we formulate and solve the robust OCP with fixed feedback $K=0$ as a standard OCP with augmented state.
This serves as a loose lower bound on the computational time for the more challenging OCP with feedback optimization.
For each variant, we perform five closed-loop MPC simulations of 40 time steps each, and average the resulting computation times.
Each MPC solver is fully converged for every time step.
The results are shown in Fig.~\ref{fig:timings_order_nx}.
The Riccati recursion adds very little computational cost compared with the nominal OCP and ZORO.
Due to the added nonlinearity from the trajectory-dependent uncertainty weighting, the constraint-adaptive Riccati-ZORO variant needs slightly more iterations than ZORO and Riccati-ZORO with a constant Hessian.

\subsection{Collision avoidance (with \texttt{acados})}
As a practical case study, we consider collision avoidance MPC for a differential drive robot.
This problem formulation has previously been employed for the real-world application of ZORO-MPC in \cite{Gao2023}.
For this problem, SIRO exhibits unreliable convergence depending on the specific instances of \eqref{ocp:rob} encountered throughout the MPC loop.

The dynamic model uses the state $x=(p_\mathrm{x}, p_\mathrm{y}, \theta, v, \omega)$, consisting of positions $(p_\mathrm{x}, p_\mathrm{y})$, orientation $\theta$, speed $v$, and angular velocity $\omega$.
The controls are the acceleration $a$ and angular acceleration $\alpha$.
The robot should track a given reference trajectory while avoiding collision and respecting constraints on the speed, angular velocity, and actuator limits.
For more details, see \citet{Frey2024, Gao2023}.

\begin{figure}[t]
   \centering
   \includegraphics[width=\columnwidth]{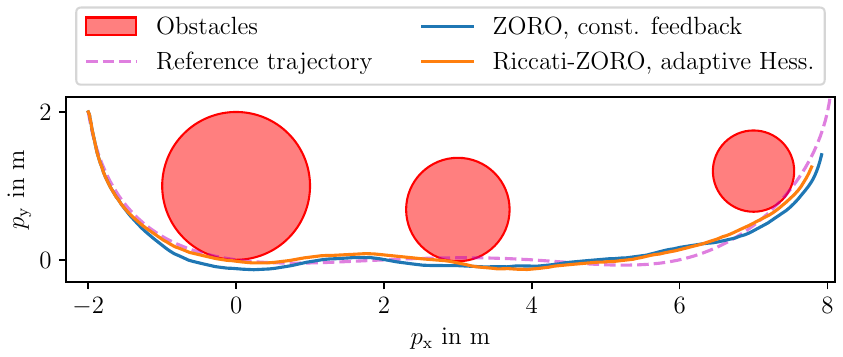}\\[5pt]
   \includegraphics[width=\columnwidth]{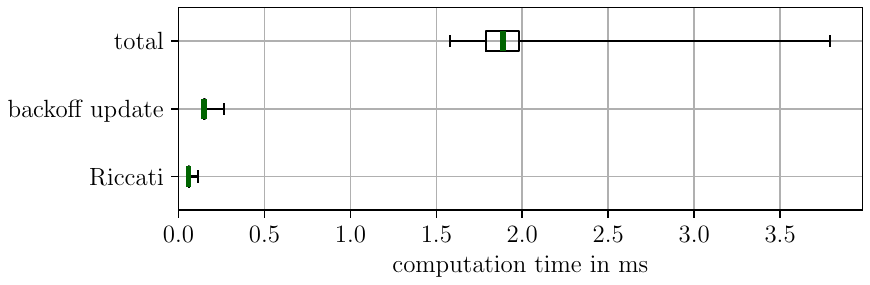}
   \caption{Top: Closed-loop trajectories for the differential drive collision avoidance case study. The online feedback optimization of Riccati-ZORO enables closer tracking of the reference in the vicinity of obstacles, resulting in a minimum obstacle distance of 0.009~m, compared to 0.058~m for ZORO with precomputed feedback.
   Bottom: Computational cost per MPC step for Riccati-ZORO.
   The boxes indicate the median and interquartile range; the whiskers the full data spread.
   The total time comprises the nominal SQP step (e.g., derivative computation and QP solve)
   and the backoff update.
   The backoff update comprises the Riccati recursion, ellipsoid propagation, and backoff computation from the ellipsoids.
   The added cost of the backoff update in comparison with the nominal SQP step is only a fraction of the latter.%
   }
   \label{fig:diff-drive-all}
\end{figure}

We compare two MPC policies: (i) ZORO with a constant feedback gain, which has been designed specifically for the differential drive robot \citep{Gao2023}, (ii) Riccati-ZORO with an adaptive Hessian.
For each MPC policy we use the real-time iteration (RTI) variant \citep{Diehl2005}, with two SQP iterations per time step, and a prediction horizon of $N=20$.
The results are shown in Fig.~\ref{fig:diff-drive-all}.

The online optimization of the feedback gains enables improved reference tracking in the vicinity of obstacles, while the computational cost is only marginally increased by the Riccati recursion.
Furthermore, note that some design effort went into the application-specific feedback gain design, whereas the Riccati-ZORO approach is generic.
In general, for nonlinear systems no predesigned feedback gain is available.

\section{Conclusion}
\label{sec:conclusion}

We have presented the algorithm Riccati-ZORO for heuristic feedback gain optimization in tube-based robust and stochastic MPC.
Riccati-ZORO reliably converges to feasible solutions of robustified MPC problems with computational costs comparable to those of nominal MPC.
Although the key principle is relevant to a variety of nonlinear tube OCPs, we have developed it specifically for ellipsoidal tubes with linearization-based dynamics.
Future work should extend it to a broader range of formulations.
Of specific interest are tube OCPs that provide robustness guarantees for the nonlinear case, and those using a learning-based uncertainty quantification.
Furthermore, due to its close connection to SIRO, the heuristic multiplier smoothing could form the basis for exact feedback gain optimization with broader convergence guarantees.

\printcredits

\bibliographystyle{cas-model2-names}

\bibliography{syscop}

\end{document}